\documentclass[12pt,leqno]{amsart}  
\usepackage{amsmath,amstext,amsthm,amssymb,amsxtra}
\usepackage[colorlinks,citecolor=red,pagebackref,hypertexnames=false]{hyperref}
\usepackage[backrefs]{amsrefs}
\allowdisplaybreaks
\swapnumbers
 
\theoremstyle{plain} 
 
\newtheorem{proposition}[equation]{Proposition} 
\newtheorem{theorem}[equation]{Theorem}

\theoremstyle{definition}

\theoremstyle{remark}
\newtheorem{remark}[equation]{Remark}

\numberwithin{equation}{section}

\def\norm#1.#2.{\lVert#1\rVert_{#2}}
\def\Norm#1.#2.{\bigl\lVert#1\bigr\rVert_{#2}}
\def\NOrm#1.#2.{\Bigl\lVert#1\Bigr\rVert_{#2}}
\def\NORm#1.#2.{\biggl\lVert#1\biggr\rVert_{#2}}
\def\NORM#1.#2.{\Biggl\lVert#1\Biggr\rVert_{#2}}

\def\ip#1,#2,{\langle #1,#2\rangle}
\def\Ip#1,#2,{\bigl\langle#1,#2\bigr\rangle}
\def\IP#1,#2,{\Bigl\langle#1,#2\Bigr\rangle}

\def\mid{\,:\,}

\def\abs#1{\lvert#1\rvert}

\def\XXint#1#2#3{{\setbox0=\hbox{$#1{#2#3}{\int}$}
     \vcenter{\hbox{$#2#3$}}\kern-.5\wd0}}

\def\eqdef{\stackrel{\mathrm{def}}{{}={}}}

\DeclareFontFamily{U}{wncy}{}
\DeclareFontShape{U}{wncy}{m}{n}{<->wncyr10}{}
\DeclareSymbolFont{mcy}{U}{wncy}{m}{n}
\DeclareMathSymbol{\Sh}{\mathord}{mcy}{"58}

\begin{document}
\title {Iterated Riesz Commutators: a simple proof of boundedness}
\author{Michael T. Lacey\and Stefanie Petermichl \and Jill C. Pipher \and Brett D. Wick}

\address{Michael T. Lacey\\
School of Mathematics\\
Georgia Institute of Technology\\
Atlanta, GA 30332\\
USA
}
\email{lacey@math.gatech.edu}
\thanks{Research supported in part by a National Science Foundation Grant.}

\address{Stefanie Petermichl\\ D\'epartement de Math\'ematiques\\ Universit\'e de Bordeaux 1\\ 351 Cours de la Lib\'eration\\ 33405 Talence, France}
\email{stefanie@math.u-bordeaux1.fr}

\address{Jill C. Pipher\\ Department of Mathematics\\ Brown University\\ Providence, RI 02912\\ USA}
\email{jpipher@math.brown.edu}
\thanks{Research supported in part by a National Science Foundation Grant.}

\address{Brett D. Wick\\ Department of Mathematics\\ University of South Carolina\\ Columbia, SC 29208\\ USA}
\email{wick@math.sc.edu}
\thanks{Research supported in part by a National Science Foundation.}

\begin{abstract}
We give a simple proof of $L^p$ boundedness of iterated commutators of Riesz transforms and a product BMO function. We use a representation of the Riesz transforms by means of simple dyadic operators - dyadic shifts - which in turn reduces the estimate quickly to paraproduct estimates.
\end{abstract}

\maketitle

\section{Introduction}

It was shown by the authors in \cite{LPPW} that product BMO of S.-Y. A Chang and R. Fefferman, 
defined on the space $\mathbb R^{d_1}\otimes \cdots 
\otimes \mathbb R ^{d_t}$, 
can be characterized by the iterated commutators of Riesz transforms
in the multi-parameter setting. This extended a classical 
one-parameter  result of R.~Coifman, R.~Rochberg, and G.~Weiss  
\cite{CRW}, and at the same time extended the work of 
M.~Lacey and S.~Ferguson \cite{sarahlacey} and M.~Lacey and E.~Terwilleger
\cite{math.CA/0310348}, which treats the case of iterated commutators 
with Hilbert transforms.  

In the multi parameter setting, one seeks to prove an inequality of the type
$$\|b\|_{BMO} \stackrel{(1)}{\lesssim}\|[T_1,[T_2,...[T_n,M_b]...]]\|_{L^2\to L^2} \stackrel{(2)}{\lesssim}\|b\|_{BMO}$$
where $T_i$ are either Riesz or Hilbert transforms acting in the $i$th variable and $M_b$ is multiplication by $b$. We refer to (1) as a lower bound and (2) as an upper bound.

In the case of the Hilbert transforms, the upper bound is not hard to prove, by appealing to analyticity 
even in the multi parameter setting,
see \cite{fergusonsadosky}, whereas 
the lower bound is extremely difficult and makes use of the corresponding upper
bound. 
Similar considerations for the Riesz transforms fail because of a lack of analyticity. The upper bound is already 
difficult, even  in the case of one parameter, and the proof of the lower bound is non-trivial. 
Coifman, Rochberg and Weiss \cite{CRW} used 
a sharp function method to prove their one-parameter analog of Theorem~\ref{upper}.  
This method, versatile as it is in the one-parameter case, 
has  no clear analog in the product setting;  see 
Lacey \cite{math.CA/0502336} for a discussion of this in a similar context. 
The methods  of \cite{CRW} do not extend  to the multi parameter setting. 

In \cite{LPPW} we find it necessary to prove an upper bound for more generalized smooth cone projection operators to simulate the analyticity we do not have. Such operators, even in the one parameter case, are not included in the work in \cite{CRW}. The upper bound for the Riesz transforms in particular allow for a relatively simple proof that extends naturally to the multi parameter setting. In this note, we present this proof. Here is the theorem:

\begin{theorem}\label{upper}
$$\|[R_{1,j_1},[R_{2,j_2},...[R_{n,j_n},M_b]...]]\|_{L^2\to L^2}\lesssim \|b\|_{BMO}$$
\end{theorem}

Here $R _{s,\, j}$ denotes the $ j$th Riesz transform acting on the $s$th variable on $ \mathbb R ^{d_s}$ and $M_b$ denotes multiplication by the symbol $b$. The BMO norm used is Cheng Fefferman product BMO.

The proof given here is rather different from that of Coifman, Rochberg 
and Weiss.  
The proof of the upper bound presented here
follows from the decomposition of the commutators into a sum of simpler terms. 
The decomposition we allude to is naturally described in terms of 
`Haar shifts,' a method of proof used by Petermichl \cite{P1}
to address a subtle question concerning the Hilbert transform in a context with matrix valued weights and then extended by Petermichl, Treil and Volberg in \cite{PTV}  to the  Riesz transforms. This method reduces the main issue to paraproduct estimates quickly.
The paraproducts that arise are of multi-parameter form.  The specific 
result needed is Theorem~\ref{t.paraproducts} below. 
This result is due to Journ\'e \cite{MR88d:42028}; 
more recent discussions of paraproducts  are in  \cites{camil1,camil2,math.CA/0502334}.

 \section{The Upper Bound }

 Let $ \operatorname M _b \varphi \eqdef b \varphi $ be the operator of pointwise multiplication 
 by a function $ b$.  
 For Schwartz functions $ f$ on  $ \mathbb R^{d}$, let $ \operatorname R_j f$ denote the 
$ j$th Riesz transform of $ f$, for $ 1\le j\le d$. 

We are concerned with product spaces $ \mathbb R ^{\vec d}=\mathbb R^{d_1}\otimes \cdots \otimes \mathbb R ^{d_t}$ for vectors $ \vec d=(d_1,\ldots , d_t)\in
\mathbb N ^{t}$. 
For Schwartz functions $b,f$ on $ \mathbb R ^{\vec d}$, 
and for a vector $\vec j=(j_1,\ldots, j_t)$ with $1\leq j_s\leq d_s$ for $s=1,\ldots,
t$ we consider the family of commutators
\begin{equation} \label{e.commutators}
\operatorname C _{\vec j}(b, f)(x)
\eqdef [\cdots[[\operatorname M_b, \operatorname R _{1,\, j_1}], 
\operatorname R _{2,\, j_2}],\cdots ](f)(x)
\end{equation}
where  $ \operatorname R _{s,\, j}$ denotes the $ j$th Riesz transform
acting on $ \mathbb R ^{d_s}$.

\begin{theorem}\label{t.upper}  We have the estimates below, valid for $ 1<p<\infty $. 
\begin{equation}\label{e.BMO}
\norm \operatorname C_{\vec j}(b,\varphi ).p. \lesssim  \norm b.\textup{BMO}.\norm \varphi .p.\,.
\end{equation}
By $ \textup{BMO}$, we mean Chang--Fefferman $ \textup{BMO}$. 
\end{theorem}

\section{Haar Functions in Several Dimensions} 

We start in the one-parameter setting. 
We will use dilation and translation operators on $ \mathbb R ^{d}$
\begin{align}\label{e.trans}
\operatorname {Tr} _{y} f(x) &\eqdef f(x-y)\,,\quad y\in \mathbb R ^{d}, 
\\
\label{e.dil}
\operatorname {Dil} _{a} ^{(p)} f(x) & \eqdef  a ^{-d/p} f(x/a)\,, 
\quad a>0\,,\ 0<p\le \infty \,.
\end{align}
These will also be applied to sets, in an obvious fashion, in the case of $ p= \infty $.  

By the \emph{canonical dyadic grid} in $ \mathbb R ^{d}$ we mean the collection of cubes 
\begin{equation*}
\mathcal D _{\text{cncl}} \eqdef \bigl\{ j 2 ^{k}+[0,2^k)^ d\mid j \in \mathbb Z ^{d}\,,\ 
k\in \mathbb Z \bigr\} 
\end{equation*}
By a \emph{dyadic grid} we mean any of the collections 
\begin{equation*}
\mathcal D_d=\mathcal D^{t,y}_d = \bigl\{ \operatorname {Dil} _{t} ^{(\infty )} 
\operatorname {Tr} _{y} I \mid I \in \mathcal D _{\text{cncl}} \bigr\}\,,
\qquad 1\le t \le 2\,,\ y\in \mathbb R ^{d}\,. 
\end{equation*}
Note that $\mathcal D _{\text{cncl}} = \mathcal D^{1,0}_d$. 
If we emphasize the role of dimension, we refer to these collections as 
\emph{$ d$-dimensional dyadic grids. }

Haar functions on $ \mathbb R ^{d}$ are now  described. For $ \varepsilon \in \{0,1\}$, set 
\begin{equation*}
h ^{0} \eqdef   -\mathbf 1 _{[-\tfrac12,0)}+\mathbf 1 _{[0,\tfrac12)}\,, 
\qquad 
h ^{1} \eqdef \mathbf 1 _{(-\tfrac12,\tfrac12)}\,.
\end{equation*}
Here, we put the superscript $ {}^0$ to denote that `the function has mean $ 0$,' 
while a superscript $ {}^1$ denotes that `the function is an  $ L^2$ normalized indicator function.'
In one dimension, for an interval $ I$, set 
\begin{equation*}
h ^{\varepsilon } _{I} \eqdef \operatorname {Tr} _{c(I)} \operatorname {Dil} _{\abs{ I}} 
^{(2)} h ^{\varepsilon }\,.
\end{equation*}
Of course for any choice of one dimensional grid $ \mathcal D$, the collections 
of functions $ \{ h ^{0} _{I}\mid I\in \mathcal D\}$  form a Haar basis for $ L^p (\mathbb R )$.

Let $\textup {Sig}_{d} \eqdef \{0,1\}^d-\{\vec 1\}$, which we refer to as 
\emph{signatures}.
In $ d$ dimensions, for a cube $Q$ with side $I$, i.e.,  $Q=I\times \cdots \times I$, and a choice of $ \varepsilon\in\textup{Sig}_d$, set 
\begin{equation*}
h ^{\varepsilon} _Q (x_1,\dotsc,x_d)\eqdef \prod _{j=1} ^{d} h _{I} ^{\varepsilon _j} (x_j).
\end{equation*}
It is then the case that the collection of functions 
\begin{equation*}
\operatorname {Haar} _{\mathcal D_d} \eqdef \{ h _{Q} ^{\varepsilon } 
\mid Q\in \mathcal D_d\,, \ \varepsilon\in\textup{Sig}_d\}
\end{equation*}
form a Haar basis for $ L^p (\mathbb R ^{d})$ for any choice of  $ d$-dimensional 
dyadic grid $ \mathcal D_d$. 
Here, we are using the notation $ \vec 1=(1,\dotsc,1)$.  While we exclude the 
superscript $ {} ^{\vec 1}$ here, they play a role in the theory of 
paraproducts.

We will use these bases in the tensor product setting. Thus, for 
a vector $ \vec d=(d_1,\dotsc,d_t)$, and $ 1\le s\le t$, 
let $ \mathcal D _{d_s}$ be a choice of $ d_s$ dimensional dyadic grid, and let 
\begin{equation*} 
\mathcal D _{\vec d}=\otimes _{s=1} ^{t} \mathcal D _{d_s}\,. 
\end{equation*}
Also, let $\textup {Sig}_{\vec d} \eqdef \{\vec\varepsilon=(\varepsilon_1,\ldots,\varepsilon_t):\varepsilon_s\in\textup{Sig}_{d_s}\}$.
Note that each $\varepsilon_s$ is a vector, and so $\vec\varepsilon$ is a `vector of vectors'.  For a rectangle $ R=Q_1\times \cdots Q_t$, i.e. a product of cubes 
of possibly different dimensions, and a choice of vectors $ \vec\varepsilon\in\textup {Sig}_{\vec d}$
set 
\begin{equation*}
h _{R} ^{\vec \varepsilon }(x_1,\dotsc,x_t)=\prod _{s=1}^t h _{Q_s} ^{\varepsilon _s}(x_s)
\end{equation*}
These are the appropriate functions and bases to analyze multi-parameter 
paraproducts and commutators.  

Let
\begin{equation*}
\operatorname {Haar} _{\mathcal D _{\vec d}} \eqdef \bigl\{ h _{R} ^{\vec \varepsilon }\mid 
R\in \mathcal D _{\vec d}\,, \ \vec\varepsilon\in\textup {Sig}_{\vec d}\bigr\}\,. 
\end{equation*}
This is a basis in $ L^p(\mathbb R ^{\vec d})$, where we use the notation 
\begin{equation*}
\mathbb R ^{\vec d} \eqdef \mathbb R ^{d_1}\otimes \cdots \otimes \mathbb R ^{d_t}
\end{equation*}
to emphasize that we are in a tensor product setting. 

\section{Chang--Fefferman $ \textup{BMO}$} \label{s.cfBMO}

We describe the elements of product Hardy space theory, 
as developed by S.-Y.~Chang and R.~Fefferman \cites{cf1,cf2,
MR90e:42030,MR86f:32004,MR81c:32016} as well as Journ\'e \cites{MR87g:42028,MR88d:42028}.  
By this, we mean 
the Hardy spaces associated with domains like $\otimes_{s=1}^t \mathbb R^{d_s}$.

\begin{remark}\label{r.H(Rd)}
The (real) Hardy space $ H^1 (\mathbb R ^{d})$ typically denotes the class of functions 
with the norm 
\begin{equation*}
\sum _{j=0} ^{d} \norm \operatorname R_j f.1.
\end{equation*}
where $ \operatorname R_j$ denotes the $ j$th Riesz transform. 
Here and below we adopt the convention that $ \operatorname R_0$, 
the $ 0$th Riesz transform, is the identity. 
This space is 
invariant under the one-parameter family of isotropic dilations, while   $ H^1(\mathbb R^{\vec d}) $ 
is invariant under dilations of each coordinate separately. That is, it is invariant 
under a $t$-parameter family of dilations, hence the terminology  `multiparameter' 
theory.
\end{remark}

As before, the space $H^1(\mathbb R^{\vec d}) $ has a variety of equivalent norms, 
in terms of square functions, maximal functions and Riesz transforms.  For our discussion, the characterization in terms of Riesz transforms is most useful:

$$
\norm f. H^1 (\mathbb R ^{\vec d}). =
\sum_{\substack{\vec 0\le \vec \jmath \le \vec d }}   
\NOrm \prod _{s=1} ^{t} \operatorname R_{s,j_s} f.1.\,.
$$

 $ \operatorname R_{s,j_s}$ is  the Riesz transform 
computed in the $j_s$th direction of the $s$th variable, and the 
$ 0$th Riesz transform is the identity operator.

\subsection{$ \textup{BMO} (\mathbb R^{\vec d})$}

The dual of the real Hardy space is $ H^1(\mathbb R^{\vec d})^\ast=\text{BMO}(\mathbb R^{\vec d})$,
the $t$--fold product $\text{BMO}$ space. It is a Theorem of S.-Y.~Chang and R.~Fefferman 
\cite{cf2} that this space has a characterization 
in terms of a product Carleson measure.

Define 
\begin{equation} \label{e.BMOdef}
\lVert b\rVert_{\text{BMO} ( \mathbb R^{\vec d})}\eqdef
\sup_{ U\subset \mathbb R^{\vec d}} \Bigl[\abs{   U}^{-1} 
\sum_{R\subset U}\sum_{\vec{\varepsilon}\in\textup {Sig}_{\vec d}}\abs{\ip b,w_R^{\vec\varepsilon},}^2\Bigr] ^{1/2} .
\end{equation}
Here the supremum is taken over all open subsets $U\subset \mathbb R^{\vec d}$ with finite measure, 
and we use a wavelet basis $ w_R ^{\vec \varepsilon }$.

\begin{theorem}[Chang--Fefferman $ \textup{BMO}$]\label{t.changfefferman}
We have the equivalence of norms 
\begin{equation*}
\norm f. (H^1 (\mathbb R ^{\vec d})) ^{\ast}. 
\approx \norm f. \textup{BMO}(\mathbb R ^{\vec d}). .
\end{equation*}
That is, $ \textup{BMO}(\mathbb R ^{\vec d})$
is the dual to $H^1 (\mathbb R ^{\vec d})$.
\end{theorem}

\section{Paraproducts} 

We recall a result of Journ\'e \cite{MR88d:42028}, with a recent 
improvement of Muscalu, Pipher, Tao and Thiele \cites{camil1,camil2}. (Also see 
Lacey and Metcalfe \cite{math.CA/0502334}.)

Consider the bilinear operators, in fact  multi-parameter paraproducts, on functions 
in $ \mathbb R ^{\vec d}$, given as 
\begin{equation*}
\operatorname B(f_1,f_2) \eqdef \sum _{R\in \mathcal D _{\vec d}}  \epsilon _{R}
\ip f_1, h _{R} ^{\varepsilon _1}, \ip f_2, h _{R} ^{\varepsilon _2}, 
\frac{ h _{R}  ^{\varepsilon _3} } {\sqrt{ \abs{ R} } }\,, \qquad \epsilon _R\in\{-1,+1\}
\end{equation*}
we suppress the dependence of this operator on the three indices $ \varepsilon _j\in 
\{0,1\}^{\vec d}$, as well as the choices of signs $ \epsilon _R$. 

\begin{theorem}\label{t.paraproducts} 
Recall that $ \vec d=(d_1,\dotsc,d_t)$ and that $ \varepsilon _j=(\varepsilon _{j,1},\dotsc, 
\varepsilon _{j,t})$.  If for all $ 1\le s\le t$, there 
is at most one  choice of $ j=1,2,3$ with $ \varepsilon _{j,s}=\vec1$, 
then the operator $ \operatorname B$ satisfies 
\begin{equation*}
\operatorname B\mid L^p \times L ^{q}
\longrightarrow L^r\,, 
\qquad 1<p,q<\infty \,, \ \tfrac1p+\tfrac1q=\tfrac1r\,.
\end{equation*}
If in addition, $ \varepsilon _1\neq\vec 1$, we will have the estimates 
\begin{equation*}
\operatorname B\mid\textup{BMO}û\times L ^{p}
\longrightarrow L^p\,, \qquad 1<p<\infty \,. 
\end{equation*}
\end{theorem}

\section{Discrete Analogs of Riesz Transforms} 

The Riesz transforms can be obtained from certain types of linear operators 
which map Haar functions to themselves (Haar transforms), a fact demonstrated by 
 Petermichl \cite{P1}  for the Hilbert transform, 
 and Petermichl, Treil and Volberg 
\cite{PTV} for Riesz transforms.

Let us describe the Haar transforms we consider. 
 Fix a dimension $ d$, and a choice of $ d$-dimensional
dyadic  grid $ \mathcal D$.  Let $\sigma \mid \mathcal D\rightarrow \mathcal D $ such 
that $ 2^d \abs{ \sigma (I)}=\abs{ I}$, for all $ I\in \mathcal D$.  Use the same 
notation for a map 
\begin{equation*} 
\sigma \mid \textup{Sig}\longrightarrow \textup{Sig}\cup\{0\}
\end{equation*}
where $ \textup {Sig} \eqdef \{0,1\}^d-\{\vec 1\}$ is the set of `signatures' 
associated with the wavelets. If $\sigma(\epsilon)=0 $ then $h^{\sigma(\epsilon)} :=0$. 

Define $ \operatorname Q $ by 
\begin{equation}\label{e.Qdef} 
\operatorname Q h ^{\varepsilon }_I \eqdef h ^{\sigma (\varepsilon )} _{\sigma (I)},
\end{equation}
and then extending linearly.  
There are two facts that we need about these transforms.  See \cite{PTV} for the proof of this proposition.

\begin{proposition}\label{p.Q}
(a) The operators $ \operatorname Q$ as defined above 
map $ L^p(\mathbb R ^{d})$ into itself for all $ 1<p<\infty $. 
\\
(b)  The Riesz transforms are in the convex hull of the class of  operators $ \operatorname Q$, the 
convex hull taken with respect to the strong operator topology. 
\end{proposition}

We need tensor products of the operators $ \operatorname Q$ we have just described. 
For $ \vec d=(d_1,d_2,\dotsc,d_t)$, and $ 1\le s\le t$, let $ \operatorname Q_s$ 
be an operator as above, acting on $ L^2(\mathbb R ^{d_s}) $.  We will use the 
same notation for the operator on $ L^2(\mathbb R ^{\vec d})$ that equals 
$ \operatorname Q_s$ on the $ L^2(\mathbb R ^{d_s}) $, and is the identity on the 
orthocomplement of $ L^2(\mathbb R ^{d_s}) $ in $ L^2(\mathbb R ^{\vec d})$. 

\begin{proposition}\label{p.Qvec} The operators 
\begin{equation} \label{e.Qvec}
\vec{\operatorname Q} \eqdef \operatorname Q_1\otimes \cdots \otimes\operatorname Q_t
\end{equation}
extend to  bounded linear operators from $ L^p(\mathbb R ^{\vec d})$ to itself for 
all $ 1<p<\infty $. 
\end{proposition}

\begin{proof}[Proof of Proposition~\ref{p.Qvec}]
The proof follows by duality.  Suppose $1<p<\infty$ and $f\in L^p(\mathbb R^{\vec d})\cap L^2(\mathbb R^{\vec d})$ and $g\in L^{p'}(\mathbb R^{\vec d})$.  Computation gives
\begin{eqnarray*}
\abs{\ip \vec{\operatorname Q}f, g,} & = & \sum_{\vec{\varepsilon}\in\textup{Sig}_{\vec d}} \sum _{R\in \mathcal D _{\vec d}}  
\ip f, h^{\vec{\varepsilon}}_R, \ip g ,h ^{\sigma(\vec{\varepsilon})} _{\sigma (R)},\\
& \leq & \int_{\mathbb R^{\vec d}}{\operatorname S}_{\vec{d}}(f)(x) {\operatorname S}_{\vec{d}}(g)(x)dx
\end{eqnarray*}
where we have defined the square function
$$
{\operatorname S}_{\vec{d}}(f)(x)\eqdef\left(
\sum_{\vec{\varepsilon}\in\textup{Sig}_{\vec d}} 
\sum _{R\in \mathcal D _{\vec d}}\abs{\ip f,h^{\vec{\varepsilon}}_R, } ^2 \frac{1_R(x)}{\abs{R}}\right)^{1/2}.
$$
Using that 
\begin{equation*}
\operatorname {Haar} _{\mathcal D _{\vec d}} \eqdef \bigl\{ h _{R} ^{\vec \varepsilon }\mid 
R\in \mathcal D _{\vec d}\,, \ \vec\varepsilon\in\textup{Sig}_{\vec d}\bigr\}\,. 
\end{equation*}
is a basis for $L^p(\mathbb R^{\vec d})$ it is straightforward to show that $\norm {\operatorname S}_{\vec{d}}(f).p.$ is an equivalent norm for $L^p(\mathbb R^{\vec d})$ when $1<p<\infty$.  This and an application of Cauchy-Schwarz implies
$$
\abs{\ip \vec{\operatorname Q}f, g,} \lesssim \norm f .p. \norm g.p'..
$$
The above inequality and the density of $L^2(\mathbb R^{\vec d})\cap L^p(\mathbb R^{\vec d})$ in $L^p(\mathbb R^{\vec d})$ gives that $\vec{\operatorname Q}$ is a bounded operator.
\end{proof}

\section{Proof of Upper Bound for Commutators} 

To prove our main Theorem, it suffices to consider the  commutators 
in which the Riesz transforms are replaced by choices of the Haar transforms 
as given in (\ref{e.Qdef}). 

Define commutators by 
\begin{equation}\label{e.CQ}
\operatorname C_{\vec{\operatorname Q}}(b,f) 
\eqdef [\cdots[\operatorname M _{b},\operatorname Q_1],\cdots,\operatorname Q_t]
\end{equation}
where the $\operatorname Q_s $ are given as in (\ref{e.Qvec}), and we view the 
commutator above as acting on $ L^p( \mathbb R ^{\vec d})$. 
In the remainder of this section we prove this proposition. 

\begin{proposition}\label{p.CQ} 
The commutators $ \operatorname  C_{\vec{\operatorname Q}}$ map $ \textup{BMO}\times L^p$ into 
$ L^p$ for $ 1<p<\infty $. 
\end{proposition}

We make some detailed remarks about the one parameter case. Let $ \operatorname Q$ be as 
in (\ref{e.Qdef}), and 
consider the commutator 
\begin{equation} \label{e.oneParaComm}
[ \operatorname M _{h _{I'} ^{\varepsilon '}}, \operatorname Q] h _{I} ^{\varepsilon}
=h _{I'} ^{\varepsilon '}h ^{\sigma (\varepsilon )} _{\sigma (I)}- 
	\operatorname Q h _{I'} ^{\varepsilon '}h ^{\varepsilon } _{I}
\end{equation}
There is no contribution if $I\cap I' = \emptyset$.  The  other cases yield what follows:
\begin{equation} \label{e.cases}
[ \operatorname M _{h _{I'} ^{\varepsilon '}}, \operatorname Q]h _{I} ^{\varepsilon}
=
\begin{cases}
0  &   I\subsetneq I' 
\\
\pm\abs{I} ^{-1/2}  h ^{\sigma (\varepsilon)} _{\sigma (I)} 
-\operatorname Q (h ^{\epsilon'} _{I} h ^{\epsilon} _{I} )  
& I=I' 
\\
h ^{\varepsilon'} _{\sigma (I)}h ^{\sigma(\varepsilon)} _{\sigma(I)}
\pm \abs{I} ^{-1/2} h ^{\sigma (\varepsilon ')} _{\sigma ^2(I)} 
&  I'=\sigma (I) 
\\
\pm \abs{I} ^{-1/2} h ^{\sigma (\varepsilon ')} _{\sigma (I')}
& I'\subsetneq I \textup{ and } I' \cap \sigma(I)=\emptyset
\\   
\pm\abs{\sigma(I)} ^{-1/2}h ^{\varepsilon '} _{I'}
\pm \abs{I}^{-1/2} h ^{\sigma (\varepsilon ')} _{\sigma (I')} 
& I'\subsetneq \sigma(I)\,.
\end{cases}
\end{equation}

The first triangular sum corresponding to $I\subsetneq I'$ is therefore trivial. It illustrates the essential cancellation contained in the commutator. It remains to consider the diagonal sums $I=I'$ and $\sigma(I)=I'$ as well as the other triangular sum $I'\subsetneq \sigma(I)$. All these do not require any cancellation of the commutator as we shall see.

First, we calculate the diagonal part for $I=I'$ explicitly and obtain the two sums:
$$\sum_I\sum_{\varepsilon, \varepsilon'\neq\vec{1}} \pm(b,h_I^{\varepsilon'})(f,h_I^{\varepsilon})|I|^{-1/2}Q(h_I^{\varepsilon}) \;\textup{ and }\;Q\left(\sum_I\sum_{\varepsilon, \varepsilon'\neq\vec{1}} (b,h_I^{\varepsilon'})(f,h_I^{\varepsilon})|I|^{-1/2}h_I^{\tilde{\varepsilon}}\right)$$
The first sum is essentially a finite combination of operators of the form $B(b,Qf)$ and hence bounded. Here we disregard the sign changes as the paraproducts are unconditionally convergent. We also note the fact that we may change the signatures of the Haar functions, which is included in our definition of paraproduct.
In the same sense, the second sum has terms of the form $Q(B(b,f))$. The term $h_I^{\tilde{\varepsilon}}$ stems from the product of two Haar functions based on $I$. The signature $\tilde{\varepsilon}$ may be equal to $\vec{1}$, in which case we still have a convergent paraproduct.

For the other diagonal sum, which corresponds to $I'=\sigma(I)$, the explicit calculation and the reasoning is similar. 

At last, we demonstrate the triangular sum by explicit calculation. We start by summing the first summand over $I$:
$$\sum_{I'}\sum_{\varepsilon'\neq\vec{1}} (b,h_{I'}^{\varepsilon'})h_{I'}^{\varepsilon'}\sum_{I:I\supsetneq I'} \sum_{\varepsilon\neq\vec{1}}
(f,h_I^{\varepsilon})Q(h_I^{\varepsilon})$$
The inner part is equal to 
$\sum_{I:I\supsetneq I'} \sum_{\varepsilon\neq\vec{1}}
(Qf,h_I^{\varepsilon})h_I^{\varepsilon}$
which on $I'$ equals a renormalized average of $Qf$, namely $(Qf, h^{\vec{1}}_{I'})|I'|^{-1/2}$. We conclude that the first sum is a combination of paraproducts of the form $B(b,Qf)$. Similarly, the second sum has terms of the form $QB(b,f)$ with an indicator falling on $f$.

This proves that our commutator is a finite linear combination of terms of the form
\begin{equation*}
\operatorname Q \operatorname B(b,f),\qquad \operatorname B(b,\operatorname Q f)
\end{equation*}
for appropriate choices of $\operatorname  Q$ and paraproducts  $\operatorname B $.  It is essential 
to note that these paraproducts are of the form that we can permit $ b\in\textup{BMO}û$, 
so that the one parameter form of Theorem~\ref{t.paraproducts} applies to prove 
Proposition~\ref{p.CQ} in this case. 
This completes the discussion in the one parameter case.

\bigskip 

In passing to the multi-parameter paraproducts, we adopt the notation of Proposition~\ref{p.Qvec}.
Thus, we have a choice of $ \vec d=(d_1,\dotsc,d_t)$ dimensional dyadic rectangles $ \mathcal D_{\vec d}$, and a choice of operators  $ \operatorname Q_s$,
for $ 1\le s\le t$, 
which is an  operator as in (\ref{e.Qdef}) acting on $ \mathbb R ^{d_s}$.  
We can then make 
an explicit computation of the commutator as in (\ref{e.oneParaComm}), namely 
\begin{equation*}
[\cdots[ \operatorname M _{h ^{\vec \varepsilon } _{R}} , \operatorname Q_1], \cdots 
\operatorname Q _t]
\end{equation*}
The result is a tensor product of the results on the right hand side of   (\ref{e.cases}). In other words, the multi parameter commutator splits into its components.

It follows that we can write the commutator as a finite linear combination of terms 
\begin{equation*}
\vec{\operatorname Q} \operatorname B(b,f)\,,\qquad \operatorname B(b,\vec{\operatorname Q}f)
\end{equation*}
for different choices of multi-parameter paraproduct $ \operatorname B $ and different choices 
of operator $ \vec{\operatorname Q}$.  
Thus Proposition~\ref{p.CQ} follows from Theorem~\ref{t.paraproducts}.

\end{document}